\documentclass[11pt]{amsart}

\usepackage[colorlinks=true,citecolor=cyan,backref=page]{hyperref}

\usepackage[square,sort&compress,comma,numbers]{natbib}
\usepackage{nicefrac,xcolor,upref,version}

\usepackage{amscd,amsthm,amsmath,amssymb,amsfonts}

\usepackage{mathrsfs}

\newtheorem{theorem}{Theorem}[section]

\newtheorem{corollary}[theorem]{Corollary}

\numberwithin{equation}{section}

\def\beq{\begin{equation}}
\def\eeq{\end{equation}}

\def\cT{{\mathcal T}}
\def\cU{{\mathcal U}}

\begin{document}

\vskip 5mm

\title[On the difference between perfect powers and integral $S$-units]
{On the difference between perfect powers and integral $S$-units}

\author{Yann Bugeaud}
\address{I.R.M.A., UMR 7501, Universit\'e de Strasbourg
et CNRS, 7 rue Ren\'e Descartes, 67084 Strasbourg Cedex, France}
\address{Institut universitaire de France}
\email{bugeaud@math.unistra.fr}

\begin{abstract}
Let $q_1, \ldots , q_t$ be distinct prime numbers. 
Let $a_1, \ldots , a_t$ be nonnegative integers. 
We establish effective lower bounds for $|z^d -  q_1^{a_1} \ldots q_t^{a_t}|$
and for its greatest prime factor, which tend to infinity with $z^d$, 
where $z$ is a positive integer 
coprime with $q_1 \ldots q_t$ and $d \ge 2$ is an integer.
\end{abstract}


\maketitle

\section{Introduction}\label{sec:1} 

It is by no means obvious that the distance between powers of $2$ and powers of $3$ 
tends to infinity. This has been proved by Thue \cite{Thue09} in 1909, but his proof does not 
allow one to determine the set of powers of $2$ which differ by a given integer from a power of $3$. 
In 1935 Gelfond \cite{Gel35} established an effective lower bound for $|3^m - 2^n|$ which tends to infinity 
as the maximum of the positive integers $m, n$ tends to infinity. Thue's result was extended by Mahler 
\cite{Mah53} who proved that the greatest prime divisor of $|3^m - 2^n|$ tends to infinity with the maximum 
of the postitive integers $m$ and $n$, but, as for Thue's result, Mahler's proof does not yield any effective 
lower bound.  The interested reader will find further references in the monographs \cite{ShTi86,Rib94}. 
All these results have been considerably generalized, and in an effective way, after the end of the 60s and the development 
of Alan Baker's theory of linear forms in an arbitrary number of logarithms of algebraic numbers (Gelfond's result 
solves the case of linear forms in two logarithms); see \cite{Bu18b} and references therein. 
In this note, we present a new extension. We replace the sequence 
of powers of  $2$ by a sequence of positive integers whose prime divisors belong to a given, finite set $T$, and the 
sequence of powers of $3$ by the sequence of perfect powers coprime with the elements of $T$.

Let $T := \{q_1, \ldots , q_t\}$ be a non-empty set of distinct prime numbers and let $\cT = (t_i)_{i \ge 1}$ denote the set of  
integers greater than $1$ whose prime divisors are all in $T$, ranged in increasing order. 
A result of Tijdeman \cite{Tij73} asserts that the difference $t_{i+1} - t_i$ tends to infinity with $i$ and, more precisely, that 
there exists an eﬀectively computable real number $C$, depending only on $T$, such that
$$
t_{i+1} - t_i \ge t_i (\log 2 t_i)^{-C}, \quad i \ge 1. 
$$
This result is essentially best possible \cite{Tij74}.

For an integer $d \ge 2$, 
let $\cU_d$  
denote the set composed of the perfect $d$-th powers $z^d$, with $z \ge 2$ 
coprime with $q_1 \ldots q_t$, ranged in increasing order. Let
$$
\cU = \bigcup_{d \ge 2} \cU_d = (u_j)_{j \ge 1}
$$
denote the union over $d \ge 2$ of the sets $\cU_d$, ranged in increasing order. 
While we know since Mih\u ailescu \cite{Mih04} that $8$ and $9$ are the only consecutive perfect powers (Catalan's conjecture), 
we cannot exclude that 
there could be infinitely many integers $j$ such that $u_{j+1} = u_j + 2$. In this note, we are concerned 
with the gap between the sets $\cT$ and $\cU$.

Let $d \ge 2$ be an integer. Let $x$ be an element of $\cT$ and $u$ an element of $\cU_d$. 
For $d \ge 3$ it follows from classical results on Thue and Thue--Mahler equations that 
$| x - u|$ and its greatest prime factor tend to infinity as the maximum of $x$ 
and $u$ tends to infinity and  effectively computable 
lower bounds have been established; see Section \ref{sec:2}. The same holds when $d=2$, as recently proved in \cite{Bu26}.
In these results, the exponent $d$ is fixed. The aim of this note is to give effective lower bounds 
independent of $d$ for $|x - u|$ and its greatest prime factor.

For an integer $z$ with $|z| \ge 2$, let $P[z]$ denote the greatest prime factor of $|z|$. 
Set $P[0] = P[-1]= P[1] = 1$. 
For any positive real number $x$, we set $\log_* x = \max\{1, \log x\}$.

\begin{theorem} \label{main} 
We keep the above notation. 
There exist effectively computable, positive, real numbers $\kappa_1, \kappa_2$, depending only on $T$, such that, 
for every $x$ in $\cT$ and every $u$ in $\cU$ we have 
$$
|x - u| \ge (\log_* \log X)^{\kappa_1} 
$$
and
$$
P[x - u] \ge \kappa_2 \,  \log_* \log_* \log X \, \frac{\log_* \log_* \log_* \log X}{\log_* \log_* \log_* \log_* \log X}, 
$$
where we have set $X := \max\{x, u\}$. 
\end{theorem}

For the proof of Theorem \ref{main}, we simply combine several results already in the literature (or somehow hidden). 

We display an immediate corollary of Theorem \ref{main}. 

\begin{corollary}  \label{cor}
We keep the above notation. 
Let $m$ be an integer. The Diophantine equation 
$$
x - u = m, \quad \hbox{in $x \in \cT, u \in \cU$},
$$
has only finitely many solutions, and all of them can be effectively determined. 
\end{corollary}

Theorem \ref{main} and its corollary illustrate the power of Baker's theory of linear forms in logarithms. 
They can be extended to arbitrary number fields without any new argument. 

Some parts of Corollary \ref{cor}, e.g. when $u$ is in $\cU$ with $d \ge 3$, follow from Theorem 3 (see also (12)) of Gy\H ory, Pink, and 
Pint\'er \cite{GyPiPi04}; see also \cite[Corollary 2.4]{GyPi08} and \cite[Theorem 12.2]{ShTi86}.

Throughout this note, the constants $C_1, C_2, \ldots$ and the constants 
implied by $\ll, \gg$ are positive, effectively computable and depend only on $T$. 
Furthermore, we let $p_1 = 2, p_2 = 3, p_3, \ldots$ denote the sequence of prime numbers ranged in increasing order 
and recall that, by the Prime Number Theorem, $p_s / (s \log s)$ tends to $1$ as $s$ tends to infinity. 

We keep the notation from Section \ref{sec:1} throughout this note. 

\section{The case of a fixed exponent}  \label{sec:2}

Let $a_1, \ldots , a_t$ be nonnegative integers, not all zero, and set 
$$
x := q_1^{a_1} \ldots q_t^{a_t}.
$$
Let $d \ge 2$ be an integer. Let $u_d$ be in $\cU_d$ 
and write $u_d = z^d$ with $z \ge 2$.

Consider first the case $d=2$. It follows from \cite[Theorem 1.4]{Bu26} and \cite[Corollary 1.3]{Bu26} that 
\beq \label{xu2}
\log |x - u_2| \gg \log X_2
\eeq
and
\beq  \label{Pxu2}
P[x - u_2] \gg \log_* \log X_2 \frac{\log_* \log_* \log X_2}{\log_* \log_* \log_* \log X_2},
\eeq
where we have set $X_2 := \max\{x, u_2\}$.

Assume now that $d \ge 3$. Write
$$
\Delta := u_d - x = z^d -  q_1^{r_1} \ldots q_t^{r_t} \, \big(q_1^{\lfloor a_1 / d \rfloor} \ldots q_t^{\lfloor a_t / d \rfloor}\big)^d,
$$
where $r_i$ is the remainder in the Euclidean division of $a_i$ by $d$, for $i=1, \ldots , t$. Thus, we get $t^d$ integer binary forms 
$$
F_{\underline r} (Z, Y) = Z^d -  q_1^{r_1} \ldots q_t^{r_t} Y^d, \quad {\underline r} = (r_1, \ldots , r_t), 
$$
of degree $d$ and $t^d$ Thue equations $F_{\underline r} (Z, Y) = \Delta$.

It then follows from the explicit bound for the solutions to Thue equations given in 
\cite[Theorem 3]{BuGy96b} that there exists $C_1$ such that  
\beq \label{xud}
\log X_d \le (C_1 d)^{18 d}  \log |\Delta|, 
\eeq
where we have set $X_d := \max\{x, u_d\}$.

We can as well derive a lower bound for the greatest prime factor of $\Delta$.
Write 
$$
u_d - x = z^d -  q_1^{r_1} \ldots q_t^{r_t} \, \big(q_1^{\lfloor a_1 / d \rfloor} \ldots q_t^{\lfloor a_t / d \rfloor}\big)^d
= \pm p_1^{b_1} \ldots p_s^{b_s},
$$
with $b_s \ge 1$, in such a way that 
$P[x - u_d] = p_s =: P$. 
Since the degree of the splitting field of the polynomial
$F_{\underline r} (X, 1)$ is bounded from above by $d (d-1)$, 
it follows from (the proof of) \cite[Theorem 4]{BuGy96b} that there are $C_2$ and $C_3$  
such that 
$$
\log X_d \le  (C_2 d s)^{14 d s} P^{d (d-1)} C_3^{2 d^2} (\log P)^{d s + 2},
$$
thus
\beq  \label{xd2P}
\log_* \log X_d \ll d^2 P. 
\eeq
Observe that the slightly better estimate 
$$
\log_* \log X_d  \ll  d^2 \frac{P \log_* \log P}{\log P} 
$$  
follows from \cite{GyYu06} (see also \cite[Theorem 2.5]{BuEvGy18})
and is of the same strength as \eqref{Pxu2}. 
In view of the bounds on $d$ obtained in Section \ref{sec:3}, this 
slight improvement on \eqref{xd2P} has no influence 
on Theorem \ref{main}.


\section{Completion of the proof of Theorem \ref{main}}   \label{sec:3}

We treat the exponent $d$ as a variable and complete the proof of Theorem \ref{main}. 

Write 
$$
q_1^{h_1} \ldots q_t^{h_t} \, \big(q_1^{\lfloor a_1 / 2 \rfloor} \ldots q_t^{\lfloor a_t / 2 \rfloor}\big)^2  + \Delta = z^d,
$$
where $h_i$ is in $\{0, 1\}$ for $i=1, \ldots , t$. 
It follows from \cite[Theorem 2.2]{BBGMOS} (which gives a fully explicit version of the 
Schinzel--Tijdeman theorem \cite{ScTi76}) that
$$
\log d \ll \log |\Delta|. 
$$ 
Combined with \eqref{xu2} and \eqref{xud}, we derive that 
$$
\log |\Delta| \gg  \log_* \log_* \log X,
$$
which proves the first assertion of Theorem \ref{main}.

Recall that
$$
\Delta := z^d - q_1^{a_1} \ldots q_t^{a_t} = \pm p_1^{b_1} \ldots p_s^{b_s},
$$
with $b_s \ge 1$ and $P = p_s$.  
We argue as in the proof of \cite[Theorem 9.1]{ShTi86}; see also \cite[Theorem 6.13]{Bu18b}. 
Set 
$$
B := \max\{a_1, \ldots , a_t, b_1, \ldots , b_s\} 
$$
and observe that 
\beq \label{dlogz}
d \log z \ll B \log P.
\eeq
Without any loss of generality, we can assume that $d$ is prime and 
that $d$ exceeds $q_1, \ldots , q_t, p_1, \ldots , p_s$. Since $d$ is prime, it is not divisible 
by $p_j$ for $j=1, \ldots , s$, nor by $q_i$ for $i = 1, \ldots , t$.  

Let $j$ be in $\{1, \ldots , s\}$ with $b_j \ge 1$. By applying 
\cite[Theorem 2.11]{Bu18b} (extracted from \cite{Yu07}) with 
$$
\alpha_n = z, \quad b_n = - d, \quad 
\{\alpha_1, \ldots , \alpha_{n-1}\} = \{q_i : a_i \ge 1, 1 \le i \le t\},
\quad p = p_j, \quad B_n = d, 
$$
and $\delta = {1 \over 2}$, we obtain 
$$
b_j \ll \, \frac{p_j}{(\log p_j)^2} \, \max \Bigl\{ (\log z)  (\log (p_j d)), {B \over d} \Bigr\}. 
$$
Arguing similarly to bound $a_1, \ldots , a_t$, we get that there exists $C_4$ such that 
$$
\max\{a_1, \ldots , a_t\} \ll \, C_4^s \, \max \Bigl\{ (\log p_1) \ldots (\log p_s) (\log z)  (\log d), {B \over d} \Bigr\}.
$$
This gives 
$$
d \ll \max \Bigl\{ \frac{p_s}{(\log p_s)^2}, C_4^s \Bigr\} 
$$
or 
$$
B \ll (\log z) (\log d) \max \Bigl\{ \frac{p_s}{\log p_s}, C_4^s (\log p_1) \ldots (\log p_s) \Bigr\}. 
$$
By \eqref{dlogz}, the latter inequality gives
$$
d \ll (\log P) (\log d) \max \Bigl\{ \frac{p_s}{\log p_s}, C_4^s (\log p_1) \ldots (\log p_s) \Bigr\}. 
$$
Consequently, we get 
$$
d \ll C_4^s \quad \hbox{or} \quad d \ll  (\log d) C_4^s (\log s)^s, 
$$
thus (recall that $P= p_s$)
$$
\log d \ll  {s \log_* \log s} \ll  P \frac{\log_* \log P}{\log P}.
$$
Combined with \eqref{Pxu2} and \eqref{xd2P}, this gives
$$
\log_* \log_* \log X  \ll {P \frac{\log_* \log P}{\log P}}, 
$$
and the proof of Theorem \ref{main} is complete.

\section*{Acknowledgement} 
The author is very thankful to K\'alm\'an Gy\H ory for several remarks on a first draft.

\end{document}